\RequirePackage{fix-cm}
\documentclass[smallextended]{svjour3}       
\smartqed  
\usepackage{graphicx}
\usepackage{color}
\usepackage{amsfonts}
\usepackage{amsmath}
\newcommand{\bu}{\textbf{\textit{u}}}
\newcommand{\bx}{{\bf x}}
\newcommand{\bia}{\textbf{\textit{a}}}
\newcommand{\bib}{\textbf{\textit{b}}}
\newcommand{\bic}{\textbf{\textit{c}}}
\newcommand{\btau}{\boldsymbol{\tau}}
\newcommand{\bphi}{{\boldsymbol \phi}}
\newcommand{\tA}{\tilde{A}}
\newcommand{\tB}{\tilde{B}}
\newcommand{\bv}{\textbf{\textit{v}}}

\begin{document}

\title{Commutation Error in Reduced Order Modeling of Fluid Flows
}


\author{B. Koc \and
	M. Mohebujjaman \and
	C. Mou \and
	T. Iliescu	
}


\institute{
	B. Koc \at
	Department of Mathematics, Virginia Tech,\\
	Blacksburg, VA, 24061, USA;\\
	\email{birgul@vt.edu}
	\and	
	M. Mohebujjaman \at
	Department of Mathematics, Virginia Tech,\\
	Blacksburg, VA, 24061, USA;\\
	\email{jaman@vt.edu}
	\and
	C. Mou \at              
	Department of Mathematics, Virginia Tech,\\
	Blacksburg, VA, 24061, USA;\\
	\email{cmou@vt.edu}	
	\and
	T. Iliescu \at 
	Department of Mathematics, Virginia Tech,\\
	Blacksburg, VA, 24061, USA;\\
	Partially supported by NSF DMS-1522656 and DMS-1821145,\\ 
	\email{iliescu@vt.edu}}           


\maketitle

\begin{abstract}
For reduced order models (ROMs) of fluid flows, we investigate theoretically and computationally whether differentiation and ROM spatial filtering commute, i.e., whether the commutation error (CE) is nonzero.
We study the CE for the Laplacian and two ROM filters: the ROM projection and the ROM differential filter.
Furthermore, when the CE is nonzero, we investigate whether it has any significant effect on ROMs that are constructed by using spatial filtering.
As numerical tests, we use the Burgers equation with viscosities $\nu=10^{-1}$ and $\nu=10^{-3}$ and a 2D flow past a circular cylinder at Reynolds numbers $Re=1$ and $Re=100$.
Our investigation shows that: 
(i) the CE exists; and 
(ii) the CE has a significant effect on ROM development for low Reynolds numbers, but not so much for higher Reynolds numbers.
\end{abstract}

\section{\uppercase{Introduction}}
\label{intro}

Reduced order models (ROMs)~\cite{benner2015survey,feppon2018dynamically,hesthaven2015certified,HLB96,noack2011reduced,quarteroni2015reduced} have been used for decades in the efficient numerical simulation of fluid flows~\cite{ballarin2016fast,ballarin2015supremizer,bergmann2018zonal,bistrian2015improved,carlberg2017galerkin,choi2017space,gunzburger2017ensemble,HLB96,majda2018model,noack2011reduced,pitton2017computational,stabile2017finite,taira2017modal,strazzullo2018model,yano2018discontinuous}.
However, when the ROM dimension is too low to capture the relevant flow features, ROMs are generally supplemented with a Correction term~\cite{baiges2015reduced,benosman2017learning,fick2017reduced,gouasmi2017priori,HLB96,kondrashov2015data,noack2011reduced,osth2014need,san2018neural,wang2012proper}.
In our recent work~\cite{xie2018data}, we have shown that this Correction term can be explicitly calculated and modeled with the available data by using the ROM projection as a spatial filter. 
We note that ROM spatial filtering has also been used to develop large eddy simulation ROMs, e.g., approximate deconvolution ROMs~\cite{xie2017approximate} and eddy viscosity ROMs~\cite{benosman2017learning,HLB96,noack2011reduced,protas2015optimal,rebollo2017certified,wang2012proper}.
In all these ROMs, it has been been assumed the differentiation and ROM spatial filtering commute:
\begin{eqnarray}
\overline{\frac{\partial u}{\partial x}}
= \frac{\partial \overline{u}}{\partial x} ,
\label{eqn:ce}
\end{eqnarray}
where $u$ is a flow variable (such as velocity) and x is a spatial direction.
In this paper, we investigate whether there exists a {\it commutation error (CE)}, i.e., whether  equality~\eqref{eqn:ce} holds.
In particular, we investigate whether there is a CE for the Laplacian, which plays a central role in fluid dynamics:
\begin{eqnarray}
\overline{\Delta u}
= \Delta \overline{u} .
\label{eqn:ce-laplacian}
\end{eqnarray}
To our knowledge, this represents the first investigation of the CE in a ROM context.

When the CE exists, we also investigate whether it has any significant effect on the ROM itself.
To this end, we consider the recently proposed {\it data-driven correction ROM (DDC-ROM)}~\cite{xie2018data}, in which the Correction term (which is generally added to improve the ROM's accuracy) is modeled using the available data~\cite{brunton2016discovering,loiseau2018constrained,lu2017data,pan2018data,peherstorfer2016data}.
To investigate the effect of the CE on the DDC-ROM, we also consider the {\it commutation error DDC-ROM (CE-DDC-ROM)}, in which available data is used to model not only the Correction term, but also the CE.
Finally, we use the {\it ideal CE-DDC-ROM (ICE-DDC-ROM)}, which is the DDC-ROM supplemented with a fine resolution representation (i.e., without any additional modeling) of the Correction term. 
When the CE-DDC-ROM and the ICE-DDC-ROM yield more accurate results than the standard DDC-ROM, we conclude that the CE has a significant effect on the DDC-ROM and, therefore, should be modeled.
As numerical tests, we use the Burgers equation with viscosities $\nu=10^{-1}$ and $\nu=10^{-3}$ and a 2D flow past a circular cylinder at Reynolds numbers $Re=1$ and $Re=100$.

The paper is organized as follows: The reduced order modeling preliminaries are provided in Section \ref{ROM-prelims}. In Section \ref{commutation-error}, a detailed derivation of the commutation error is given. In Section \ref{effect-ce-ddc}, we discuss the effects of the commutation error on ROMs that are constructed by using spatial filtering. Numerical experiments are given in Section \ref{experiments}, and conclusions and future research directions are outlined in Section \ref{conclusion}.

\section{\uppercase{Reduced Order Modeling}}\label{ROM-prelims}
To compute the ROM basis functions, we use the proper orthogonal decomposition (POD) \cite{fareed2018note,HLB96,noack2011reduced,volkwein2013proper}, which we  briefly describe in this section. We emphasize, however, that our theoretical and computational developments carry over to other ROM basis functions, such as the dynamic mode decomposition~\cite{schmid2010dynamic}. The snapshots $\{u_h^1, u_h^2, \cdots, u_h^M\}$ are the finite element (FE) solutions of at $M$ different time instances. The POD seeks a low-dimensional basis that approximates the snapshots optimally with respect to a certain norm. The commonly used $L^2$ norm will be used in this paper. The solution of the minimization problem is equivalent to the solution of the eigenvalue problem $YY^TM_h\bphi_j=\lambda_j\bphi_j$, $j=1,\cdots, N_h$, where $\bphi_j$ and $\lambda_j$ denote the vector of the FE coefficients of the POD basis functions and the POD eigenvalues, respectively, $Y$ denotes the snapshot matrix, whose columns correspond to the FE coefficients of the snapshots, $M_h$ represents the FE mass matrix, and $N_h$ is the dimension of the FE space $X^h$. The eigenvalues are real and non-negative, so they can be ordered as follows: $\lambda_1\ge\lambda_2\ge\cdots\ge\lambda_d>\lambda_{d+1}=\cdots =\lambda_{N_h}=0$, where $d$ is the rank of the snapshot matrix $Y$. The ROM basis consists of the normalized functions $\{\bphi_j\}_{j=1}^r$, which correspond to the first $r\le N_h$ largest eigenvalues. Thus, the ROM space is defined as $X^r:=\text{span}\{\bphi_1,\bphi_2, \cdots, \bphi_r\}$.

\section{\uppercase{Commutation Error (CE)}} \label{commutation-error}
As a mathematical model, we consider the incompressible time-dependent Navier-Stokes equations (NSE):
\begin{align}
\frac{\partial \bu}{\partial t}
- \nu \Delta \bu
+ \bu \cdot \nabla \bu
+ \nabla p
&= {\bf 0}\hspace{3mm} \Omega\times (0, T],
\label{eqn:nse-1}                                                         \\
\nabla \cdot \bu
&= 0 \hspace{3mm} \Omega\times (0, T],
\label{eqn:nse-2}
\\
\bu &= g(\bx) \hspace{3mm} \partial\Omega\times (0, T],        
\end{align}
where $\bu$ is the velocity, $p$ the pressure, $\nu$ the kinematic viscosity, $T$ the simulation time, and $\Omega$ the domain of the fluid. We use the initial condition $\bu(\bx, 0) = \bu_0(\bx)$. 
In this paper, we assume that $\Omega\subset\mathbb{R}^d$, $d\in\{2,3\}$, is a convex polygonal or polyhedron domain with boundary $\partial\Omega$. The discrete FE velocity and pressure spaces are denoted by $X^h$ and $Q^h$, respectively. We denote the usual $L^2(\Omega)$ norm and inner product with $\|\cdot\|$ and $(\cdot,\cdot)$, respectively.
To derive the commutation error due to filtering, we apply a continuous filter to \eqref{eqn:nse-1}. This yields the filtered-NSE (F-NSE), which have been used to develop LES models \cite{BIL05}:
\begin{align}
\displaystyle \overline{\frac{\partial \bu}{\partial t }}-\nu \overline{\Delta \bu}+ \overline{( \bu \cdot \nabla ) \bu } +\overline{ \nabla p} = {\bf 0}.\label{f-nse}
\end{align}
The F-NSE \eqref{f-nse} eliminate the small length scales in the continuous NSE \eqref{eqn:nse-1}. For this reason, a ROM for \eqref{f-nse} needs fewer POD modes than a ROM for \eqref{eqn:nse-1} to achieve a fixed numerical accuracy. However, to develop practical ROMs for the F-NSE \eqref{f-nse}, we must first investigate the commutation error, i.e., whether filtering and differentiation commute. The \textit{commutation error} (CE) for a spatial derivative is defined in \cite{BIL05} as
\begin{align}
\mathcal{E}_k[\bu](\bx) := \frac{\partial \overline{\bu}(\bx)}{\partial x_k} -\overline{\frac{\partial \bu(\bx)}{\partial x_k}},\label{comm_error_def1}
\end{align}
where $\bx$ is the spatial variable.
In this paper, we are particularly interested in the CE for the Laplacian term. Similarly to \eqref{comm_error_def1}, we define the Laplacian CE as
\begin{align}
\mathcal{E}_{\Delta}[ \bu](\mathbf{x}) := \Delta \overline{\bu(\mathbf{x})}-\overline{\Delta \bu(\mathbf{x})}.\label{comm_def2}
\end{align}

\subsection{\uppercase{ROM Spatial Filter}}\label{rom-filter}
To develop practical ROMs from the F-NSE \eqref{f-nse}, we need to replace the continuous filter in \eqref{f-nse} with discrete filters. In this paper, we use the differential and projection ROM filters. 

The \textit{ROM differential filter (DF)} \cite{wells2017evolve} is defined as: Let $\delta$ be radius of DF; for fixed $ r\le d$ and a given $\bu_d \in X^h$, the differential filter seeks $ \displaystyle \overline{\bu_d}^{DF} \in X^r $ such that
\begin{align}
\displaystyle \bigg(\Big(I-\delta^2 \Delta ~\Big)~\overline{\bu_d}^{DF},\bphi_i\bigg)=~\big(\bu_d,\bphi_i \label{DF} \big),~~~~~~~~~~~~\forall i=1,...,r.
\end{align}
By using ROM approximations for both $\overline{\bu_d}^{DF}$ and $\bu_d$ (i.e., $\overline{\bu_d}^{DF}=\sum\limits_{j=1}^r(\bia_r)_j\bphi_j$, $\bu_d=\sum\limits_{j=1}^d(\bia_d)_j\bphi_j$), we obtain the following dynamical system:
\begin{align}
\big( M_r+ \delta^2 S_r\big)~ \bia_r = M_{r\times d}~\bia_d, \label{dyn_DF}
\end{align}
where $M_r = ( \bphi_i,\bphi_j),~i,j=1,..,r$, $M_{r \times d} = ( \bphi_i,\bphi_j)~i=1,..,r,~j=1,..,d$, $S_r = ( \nabla \bphi_i, \nabla \bphi_j),~i,j=1,..,r$, and $\bia_r$ and $\bia_d$ represent coefficient vectors in $\overline{\bu_d}^{DF}$ and $\bu_d$, respectively. 

For fixed $ r\le d$ and a given $ \bu_d \in X^h$, the \textit{ROM projection filter} \cite{oberai2016approximate,wang2012proper} seeks $ \overline{\bu_d}^r \in X^r $ such that 
\begin{align}
\displaystyle \big( \overline{\bu_d}^r, \bphi_i \big)=~ \big( \bu_d, \bphi_i \big),~~~~~\forall i=1,...r. \label{PF}
\end{align}
By expanding $\overline{\bu_d}^r$ and $\bu_d$ in terms of the POD basis, we obtain the following dynamical system:
\begin{align}
M_r~ \bia_r = M_{r \times d}~ \bia_d, \label{dyn_PF}
\end{align}
where $\bia_r$ and $\bia_d$ are the coefficient vectors in $\overline{\bu_d}^r$ and $\bu_d$, respectively.

\subsection{\uppercase{Filtered-ROM}}
In Section \ref{rom-filter}, we defined two ROM spatial filters: the ROM differential filter and the ROM projection filter. In this section, we take another step in the development of practical ROMs from the F-NSE \eqref{f-nse} and replace the continuous velocity $\bu$ in \eqref{f-nse} with its \textit{most accurate approximation in the snapshot space}, i.e., with $\bu_{d}=\sum\limits_{j=1}^d(\bia_d)_j\bphi_j$, where $d$ is the rank of the snapshot matrix:
$ \forall i =1,...,r$ 
\begin{align}
\displaystyle \left( \frac{ \partial \overline{\bu_d}}{\partial t}, \bphi_i \right)- \nu ( \Delta\overline{  \bu_d} , \bphi_i ) + ( ( \overline{ \bu_d } \cdot \nabla) \overline{ \bu_d }, \bphi_i ) + \tau_i = 0, \label{fnse-weak}
\end{align}
where
\begin{align}
\displaystyle \tau_i &= -\big( \btau_d^{SFS}, \bphi_i \big),\\
\btau_d^{SFS} &=  ( \overline{ \bu_d} \cdot \nabla ) \overline{\bu_d} -\overline{ (\bu_d \cdot \nabla ) \bu_d }.\label{subfilter}
\end{align} 
The filter $-$ in the sub-filter scale stress tensor $\btau_d^{SFS}$ is either $-^r$ (i.e., the ROM projection filter \eqref{PF}) or $-^{DF}$ (i.e., the ROM differential filter \eqref{DF}).

We note that \eqref{fnse-weak}-\eqref{subfilter} is an $r$-dimensional system for the unknown $\overline{  \bu_d}\in X^r$. For clarity, we denote the unknown $\overline{  \bu_d}$ as
\begin{align}
\displaystyle \overline{\bu_d}:=\bu_r = \sum_{i=1}^{r} (\bia_r)_i \bphi_i.\label{rom-approx}
\end{align}
Using \eqref{rom-approx} in \eqref{fnse-weak}-\eqref{subfilter}, we get: $ \forall i=1,...,r $
\begin{align}
\displaystyle \left( \frac{ \partial \bu_r}{\partial t}, \bphi_i \right)+ \nu\big( \nabla  \bu_r, \nabla \bphi  \big)  + \big( ( \bu_r \cdot \nabla ) \bu_r , \bphi_i \big) +\tau_i = 0,\label{les-rom2}
\end{align}
where
\begin{align}
\displaystyle \tau_i &= -\big( \btau_r^{SFS} ,\bphi_i \big),\\
\displaystyle \btau_r^{SFS} &=  ( \bu_r \cdot \nabla) \bu_r  -\overline{( \bu_d \cdot \nabla ) \bu_d}.\label{subfilter-2}
\end{align}

Since \eqref{les-rom2} does not depend only on $\bu_r$, it is not a closed system. To close it, we need to solve the ROM \textit{closure problem}, to look for $ \overline{(\bu_d \cdot \nabla ) \bu_d }  = f(\bu_r)$. Once a ROM closure model is found, the large eddy simulation ROM (LES-ROM) \eqref{fnse-weak}-\eqref{subfilter} becomes practical. The most commonly used ROM closure models have been of eddy viscosity viscosity type \cite{wang2012proper}. Alternative ROM closure models, inspired from image processing and inverse problems (i.e., the approximate deconvolution ROM \cite{xie2017approximate}) and data-driven modeling (i.e., the data-driven correction ROM \cite{xie2018data}) have been recently proposed. We emphasize that all these LES-ROMs assume that filtering and differentiation commute. In what follows, we investigate this assumption.

\subsection{\uppercase{CE with Differential Filter}}\label{ce-dfilter}
By using the ROM differential filter, the Laplacian CE in \eqref{comm_def2} can be written as:
\begin{align}
\mathcal{E}_{\Delta}[\bu_d]:=\Delta \overline{\bu_d}^{DF} - \overline{\Delta \bu_d}^{DF},
\label{DF_error} 
\end{align}
where $-^{DF}$ represents the ROM differential filter.	We denote with
$\bia_r$, $\bia_d$, $\bib_r$ and $\bic_r$ the coefficient vectors of~ $\overline{\bu_d}^{DF}$, $\bu_d$, $ \Delta \overline{\bu_d}^{DF} $, and $ \overline{\Delta \bu_d}^{DF}$, respectively.
We start evaluating $ \Delta \overline{\bu_d}^{DF} $:
\begin{align}
\displaystyle \big( \Delta \overline{ \bu_d}^{DF},\bphi_i \big)=-\big( \nabla  \overline{\bu_d}^{DF},\nabla \bphi_i \big),~~~~~~\forall i=1,...,r.
\end{align}
Its corresponding dynamical system is
\begin{align}
\displaystyle M_r ~\bib_r = - S_r~ \bia_r. \label{prof-dy1}
\end{align}
Using \eqref{dyn_DF} and \eqref{prof-dy1} gives
\begin{align}
\displaystyle \bib_r =- M_r^{-1} S_r (M_r+\delta^2 S_r)^{-1} M_{r \times d}~ \bia_d. \label{DF_error_1}
\end{align}
Next, we evaluate $ \overline{ \Delta \bu_d}^{DF}$.
Using equation \eqref{DF}, we write the following equation for $ \overline{ \Delta \bu_d}^{DF}$:
\begin{align}
\displaystyle \bigg( \big(I-\delta^2 \Delta \big) \overline{\Delta \bu_d}^{DF}, \bphi_i \bigg)= \big( \Delta \bu_d, \bphi_i \big),~~~~\forall i=1,...,r.
\end{align}
Its corresponding dynamical system is
\begin{align}
\displaystyle \big(M_r+\delta^2S_r\big) ~\bic_r =- S_{r\times d} ~\bia_d,
\end{align}
which yields
\begin{align}
\displaystyle \bic_r = - \big(M_r+\delta^2 S_r \big)^{-1} S_{r \times d}~ \bia_d.\label{DF_error_2}
\end{align}
Since the right hand sides of equations \eqref{DF_error_1} and \eqref{DF_error_2} are not equal, we conclude that \textit{for the ROM differential filter, the Laplacian CE is nonzero}:
\begin{align}
\mathcal{E}_{\Delta}[\bu_d]  = \sum_{j=1}^{r} \big(\mathbf{b_r}-\mathbf{c_r}\big)_j \bphi_j \neq 0.
\end{align}

\subsection{\uppercase{CE with Projection Filter} }\label{ce-pfilter}
By using the ROM projection filter, the Laplacian CE in \eqref{comm_def2} can be written as
\begin{align}
\mathcal{E}_{\Delta}[\bu_d]:=\Delta \overline{\bu_d}^{r} - \overline{\Delta \bu_d}^{r} \label{PF_error},
\end{align}
where $-^r$ represents the ROM projection filter. We denote with $\bia_r$, $\bia_d$, $\bib_r$ and $\bic_r$ the coefficient vectors of the basis functions in $\overline{\bu_d}^{r}$, $\bu_d$, $ \Delta \overline{\bu_d}^{r}$, and $\overline{\Delta \bu_d}^{r}$, respectively.
We start evaluating $\Delta \overline{\bu_d}^r$ :
\begin{align}
( \Delta \overline{\bu_d}^r,\bphi_i) = -(\nabla \overline{\bu_d}^r,\nabla \bphi_i), ~~~~~~\forall i=1,...,r.
\end{align} 
Its corresponding dynamical system is
\begin{align}
M_r ~\bib_r = - S_{r }~\bia_r.\label{ce-pfdy}
\end{align}
Using \eqref{dyn_PF} and \eqref{ce-pfdy} we have
\begin{align}
\bib_r = -M_r^{-1} S_r M_r^{-1} M_{r\times d}~\bia_d \label{PF_error_1}.
\end{align}
Next, we evaluate $ \overline{\Delta \bu_d}^r $. By using (\ref{PF}) for $\overline{\Delta \bu_d}^r$, we write  the following equation for $\Delta \overline{\bu_d}^r$:
\begin{align}
(\overline{\Delta \bu_d}^r,\bphi_i) = (\Delta \bu_d ,\bphi_i) =-(\nabla \bu_d,\nabla \bphi_i), ~~~~~~\forall i=1,...,r.
\end{align}
Its corresponding dynamical system gives
\begin{align}
\bic_r = - M_r^{-1}S_{r \times d}~\bia_d, \label{PF_error_2}
\end{align}
Again, since the right hand sides of the equations \eqref{PF_error_1} and \eqref{PF_error_2} are not equal we conclude that \textit{for the ROM projection filter, the Laplacian CE is nonzero:}
\begin{equation}
\mathcal{E}_{\Delta}[\bu_d] = \sum_{j=1}^{r} \big( \bib_r-\bic_r \big)_j~\bphi_j \neq 0 . 
\end{equation}


\section{\uppercase{Effect of commutation error on DDC-ROM}}\label{effect-ce-ddc}
In this section, we investigate the effect of the commutation error on three LES-ROMs that are built from equations \eqref{les-rom2}-\eqref{subfilter-2} supplemented with the Laplacian CE \eqref{comm_def2}. The first LES-ROM that we investigate is the data-driven correction ROM (DDC-ROM) \cite{xie2018data}, which utilizes available data to construct an $r$-dimensional model for the Correction term $\btau$ in \eqref{les-rom2}-\eqref{subfilter-2}; the DDC-ROM, however, does not include a model for the Laplacian CE \eqref{comm_def2}. The second LES-ROM that we consider is the ideal CE data-driven correction ROM (ICE-DDC-ROM), which is the DDC-ROM supplemented with an exact (fine) resolution Laplacian \eqref{comm_def2} term. The third LES-ROM that we investigate is the commutation error DDC-ROM (CE-DDC-ROM), which is the DDC-ROM supplemented with an $r$-dimensional data-driven model for the Laplacian CE \eqref{comm_def2}. In Section \ref{experiments}, we investigate numerically whether the Laplacian CE has any effect on the DDC-ROM, i.e., whether the ICE-DDC-ROM and the CE-DDC-ROM yield more accurate results than the standard DDC-ROM. In this section, we outline the construction of the DDC-ROM, ICE-DDC-ROM, and CE-DDC-ROM.

First, we briefly  derive the standard Galerkin-ROM (G-ROM). The POD approximation of the velocity is defined as
\begin{align}
\bu_r(\bx,t)=\sum_{j=1}^{r}(\bia_r)_j\bphi_j(\bx),\label{grom-vel}
\end{align}
where $\{(\bia_r)_j\}_{j=1}^r$ are the sought time-depending coefficients, which are found by solving the following system of PDEs: $\forall i=1,...,r$,

\begin{align}
\displaystyle \left(\frac{\partial \bu_{r}}{\partial t},\bphi_i \right) +\nu \big(\nabla \bu_r, \nabla\bphi_i \big) + \Big( (\bu_r \cdot \nabla )\bu_r, \bphi_i \Big) = 0,\label{weak-grom}
\end{align}
where we assume that the modes $\{\bphi_1, \bphi_2, \cdots, \bphi_r\}$ are perpendicular to the discrete pressure space. This assumption holds if in snapshot creation we use standard mixed FE (such as Scott-Vogelius, or the mini-element; see, e.g., \cite{brenner2007mathematical,john2017divergence}). Plugging \eqref{grom-vel} into \eqref{weak-grom} gives the \textit{Galerkin ROM (G-ROM)}: 
\begin{align}
\dot{\bia_r} = A\bia_r + \bia_r^\top B\bia_r,\label{grom-anszt}
\end{align}
where the elements of the operators $A$ and $B$ are 
$A_{im}=-\nu(\nabla\bphi_m,\nabla\bphi_i)$ and $B_{imn}=-(\bphi_m\cdot\nabla\bphi_n,\bphi_i)$, $1\le i,~m,~n\le r$.
The CE-DDC-ROM and ICE-DDC-ROM frameworks need two steps to be constructed. In the first step, we use ROM spatial filtering to derive the \textit{exact mathematical formula} for the Correction term. To construct the DDC-ROM, ICE-DDC-ROM, and CE-DDC-ROM, we start with equations \eqref{les-rom2}-\eqref{subfilter-2}, to which we add the Laplacian CE \eqref{comm_def2}: $\forall i=1,\cdots,r$
\begin{eqnarray}
\left(\frac{\partial \bu_r}{\partial t}, \bphi_i \right) + \nu \big(\nabla \bu_r, \nabla \bphi_i \big) 
&+& \Big( (\bu_r\cdot \nabla) \bu_r, \bphi_i \Big)  
\nonumber \\
&+& \nu \big(\mathcal{E}_{\Delta}[\bu_d], \bphi_i\big) + \left(\btau_r^{SFS}, \bphi_i\right) = 0. \label{ddc-ce}
\end{eqnarray}
Equation \eqref{ddc-ce} yields the following dynamical system:
\begin{align}
\dot{\bia_r}=A\bia_r + \bia_r^\top B\bia_r+\mathcal{E}_{CE}+\btau,\label{sp-filter}
\end{align}
where $A$ and $B$ are same as in \eqref{grom-anszt} and the components of $\mathcal{E}_{CE}$ and $\btau$ are given by: $\forall i=1,...,r$
\begin{align}
(\mathcal{E}_{CE})_i &= -\nu \big( \mathcal{E}_{\Delta}[\bu_d], \bphi_i~ \big),\\
\tau_i &= - \big( \btau_r^{SFS}, \bphi_i \big).
\end{align}
To construct the DDC-ROM \cite{xie2018data}, we make the following ansatz:
\begin{align}
\btau(\bia_r)\approx\btau^{ansatz}(\bia_r)=\tilde{A}\bia_r+\bia_r^\top\tilde{B}\bia_r. \label{iddc}
\end{align}



To compute the operators $\tilde{A}$ and $\tilde{B}$ in \eqref{iddc}, we use data-driven modeling ensuring the highest accuracy of the vector $\btau$. To this end, we solve the following unconstrained optimization problem:
\begin{eqnarray}
\min_{\substack{\tA \in \mathbb{R}^{r \times r} \\[0.1cm] \tB \in \mathbb{R}^{r \times r \times r}}} \, 
\sum_{j = 1}^{M} \| \btau^{true}(t_j) - \btau^{ansatz}(t_j) \|^2
\label{eqn:DDC-ROM-least-squares-iddc} \, .
\end{eqnarray}
%


The \textit{data-driven correction ROM (DDC-ROM)} has the following form:
\begin{align}
\dot{\bia_r} = ( A+ \tilde{A})\bia_r + \bia_r^\top (B+\tilde{B})\bia_r,\label{ddc}
\end{align}
where the operators $A$ and $B$ are the G-ROM operators in \eqref{grom-anszt} and the operators $\tilde{A}$ and $\tilde{B}$ are the solution of the unconstrained minimization problem \eqref{eqn:DDC-ROM-least-squares-iddc}.

The \textit{ideal commutation error DDC-ROM (ICE-DDC-ROM)} is obtained by adding a  high-accuracy (i.e., from fine resolution numerical data) representation of the Laplacian CE \eqref{comm_def2}:
\begin{align}
\dot{\bia_r} = (A+\tilde{A})\bia_r + \bia_r^\top (B+\tilde{B})\bia_r + \mathcal{E}_{CE}.
\end{align}
To construct the CE-DDC-ROM, we make the following ansatz:
\begin{align}
(\btau + \mathcal{E}_{CE})(\bia_r) \approx (\btau + \mathcal{E}_{CE})^{ansatz}(\bia_r) = \tilde{A} \bia_r + \bia_r^\top \tilde{B} \bia_r.\label{ice-ddc-rom}
\end{align} 
To compute the operators $\tilde{A}$ and $\tilde{B}$ in \eqref{ice-ddc-rom}, we use data-driven modeling ensuring the highest accuracy of the vector $\btau + \mathcal{E}_{CE}$. To this end, we solve the following unconstrained optimization problem:
\begin{eqnarray}
\min_{\substack{\tA \in \mathbb{R}^{r \times r} \\[0.1cm] \tB \in \mathbb{R}^{r \times r \times r}}} \, 
\sum_{j = 1}^{M} \| (\btau+\mathcal{E}_{CE})^{true}(t_j) - (\btau+\mathcal{E}_{CE})^{ansatz}(t_j) \|^2
\label{eqn:DDC-ROM-least-squares-ceddc} \, .
\end{eqnarray}
The \textit{commutation error data-driven ROM (CE-DDC-ROM)} has the following form:
\begin{align}
\dot{\bia_r} = ( A+ \tilde{A})\bia_r + \bia_r^\top (B+\tilde{B})\bia_r,\label{ce-ddc}
\end{align}
where the operators $A$ and $B$ are the G-ROM operators in \eqref{grom-anszt} and the operators $\tilde{A}$ and $\tilde{B}$ are the solution of the unconstrained minimization problem \eqref{eqn:DDC-ROM-least-squares-ceddc}.

\section{\uppercase{Numerical Experiments}}\label{experiments}
In this section, we investigate numerically the following questions:\\

(Q1) Does the commutation error exist?

(Q2) If it exists, does the commutation error have a significant effect on ROMs?\\

\noindent To answer the first question, we evaluate numerically the Laplacian CE and compute its average, i.e., $\bold ||\mathcal{E}_{\Delta}[\bu_d]||_{L^2(L^2)}$,
which is calculated as follows:
$$\sqrt{\frac{{\int_{0}^{T}\|E_\Delta [\bu_d]\|^2dt}}{T}}\approx\sqrt{\frac{\sum_{j=1}^{M}(\textbf{b}_r(t_j)-\textbf{c}_r(t_j))'*M_r*( \textbf{b}_r(t_j)-\textbf{c}_r(t_j))\Delta t}{T}}.$$

To answer the second question, we test the following ROMs: the DDC-ROM \eqref{ddc}, the ICE-DDC-ROM \eqref{ice-ddc-rom}, and the CE-DDC-ROM \eqref{ce-ddc}. We emphasize that the ICE-DDC-ROM and CE-DDC-ROM include a representation of the Laplacian CE \eqref{comm_def2}, whereas the DDC-ROM does not. Thus, if the ICE-DDC-ROM and CE-DDC-ROM yield more accurate results than the DDC-ROM, we conclude that the CE plays a significant role in ROM development. In our numerical investigate, we consider two test problems: the 1D viscous Burgers equation (Section \ref{burgers}) and the 2D channel flow past a circular cylinder (Section \ref{sec:test-problem-setup}).
We compute the ROM error as the difference between the DNS solution projected onto the ROM space and the ROM solution.

\subsection{\uppercase{Experiment 1: Burgers Equation}}\label{burgers}
In our first experiment, we consider the Burgers equation: 

\begin{equation}
\begin{cases}
\displaystyle~~	u_t -\nu u_{xx} + u u_x = 0~,~~~x \in [0,1],~t\in[0,1]\\
\displaystyle~~	u(0,t)= u(1,t) = 0~,~~~t \in [0,1].
\end{cases}
\end{equation}

 The DNS results are obtained by using a linear FE scheme with mesh width $h=1/2048$ and timestep size $\Delta t=10^{-3}$. To investigate the effect of initial conditions, we consider a smooth initial condition (Section \ref{smooth-ic}) and a non-smooth initial condition (Section \ref{non-smooth-ic}). Furthermore, to investigate the effect of the viscosity parameter, we consider two viscosity values: $\nu=10^{-1}$ (Section \ref{smooth-ic} and \ref{non-smooth-ic}) and $\nu=10^{-3}$ (Section \ref{sec:burgers-lower-viscosity}).

\subsubsection{\uppercase{Smooth Initial Condition}}\label{smooth-ic}

We consider the initial condition
\begin{equation}
\displaystyle~~	u_0(x) = \frac{ 2 \nu \beta \pi sin(\pi x)}{\alpha + \beta cos(\pi x)}~,~~x \in [0,1],
\end{equation}
where $\alpha=5$ and $\beta=4$. 


\begin{table}[h!]
	\begin{center}
		\begin{tabular}{|c|c|c|c|c|c|}
			\hline
			\multicolumn{1}{ |c| }{  $r$  } &
			\multicolumn{1}{ |c| }{ $ \bold  \delta $ } &		
			\multicolumn{1}{ |c| }{ $ \bold ||\mathcal{E}_{\Delta}[u_d]||_{L^2(L^2)}   $ } & 		\multicolumn{1}{ |c| }{   $r$  } &
			\multicolumn{1}{ |c| }{ $ \bold  \delta $ } &		
			\multicolumn{1}{ |c| }{ $ \bold ||\mathcal{E}_{\Delta}[u_d]||_{L^2(L^2)}   $ } \\
			\hline
			2 &  1.00\textit{e}-01 &  1.40\textit{e}+00 & 5 &  1.00\textit{e}-01 & 1.37\textit{e}+00 \\ 
			\hline
			2 & 1.00\textit{e}-02 & 1.56\textit{e}-01 & 5 &  1.00\textit{e}-02 & 8.98\textit{e}-02 \\ 
			\hline
			2 &  1.00\textit{e}-03 & 1.53\textit{e}-01 & 5 &  1.00\textit{e}-03 & 1.26\textit{e}-03 \\ 
			\hline
			2 &  1.00\textit{e}-04 & 1.54\textit{e}-01 & 5 &  1.00\textit{e}-04 & 8.26\textit{e}-04 \\ 
			\hline
			\hline
			3 & 1.00\textit{e}-01 & 1.37\textit{e}+00 & 6 & 1.00\textit{e}-01 & 1.37\textit{e}+00 \\
			\hline 
			3 & 1.00\textit{e}-02 & 8.21\textit{e}-02 & 6 & 1.00\textit{e}-02 & 9.00\textit{e}-02 \\
			\hline
			3 & 1.00\textit{e}-03 & 3.21\textit{e}-02 & 6 & 1.00\textit{e}-03 & 9.69\textit{e}-04 \\
			\hline
			3 & 1.00\textit{e}-04 & 3.22\textit{e}-02 & 6 & 1.00\textit{e}-04 & 1.15\textit{e}-04 \\
			\hline
			\hline						
			4 &  1.00\textit{e}-01 & 1.37\textit{e}+00 & 7 & 1.00\textit{e}-01 & 1.37\textit{e}+00 \\
			\hline
			4 &  1.00\textit{e}-02 & 8.77\textit{e}-02 & 7 & 1.00\textit{e}-02 & 9.00\textit{e}-02 \\
			\hline
			4 &  1.00\textit{e}-03 &  5.54\textit{e}-03 & 7 & 1.00\textit{e}-03 & 9.63\textit{e}-04 \\
			\hline
			4 &  1.00\textit{e}-04 & 5.47\textit{e}-03 & 7 &  1.00\textit{e}-04 & 9.64\textit{e}-06 \\
			\hline
		\end{tabular}
		\caption{Burgers equation, ROM differential filter, smooth initial condition, $\nu=10^{-1}$, and $d=7$: Average CE for different $\delta$ and $r$ values.}\label{smooth-diff}
	\end{center}
\end{table}

\begin{table}[h!]
	\begin{center} 
		\scriptsize\begin{tabular}{|l|c|c|c|c|c|r|}
			\hline
			$r$ &  2 & 3 & 4 & 5 & 6 & 7 \\ \hline
			$ \bold ||\mathcal{E}_{\Delta}[u_d]||_{L^2(L^2)}   $ & 1.54\textit{e}-01 & 3.22\textit{e}-02 & 5.47\textit{e}-03 & 8.26\textit{e}-04 & 1.15\textit{e}-04 & $0$  \\ \hline		
		\end{tabular}
	\end{center}
	\caption{Burgers equation, ROM projection filter, smooth initial condition, $\nu=10^{-1}$, and $d=7$: Average CE for different $r$ values.}\label{smooth-proj}
\end{table}



First, we address (Q1), i.e., whether the CE exists.
For various $r$ and $\delta$ values, we list the CE computed by using the ROM differential filter (Table \ref{smooth-diff}) and the ROM projection filter  (Table \ref{smooth-proj}).
The main conclusion is that the CE exists for both filters, especially for low $r$ values.
We also observe that as $r$ increases, the CE decreases.
Finally, we note that, for the ROM differential filter (Table \ref{smooth-diff}), for a fixed $r$ value, as $\delta$ decreases, the CE decreases.

\begin{table}[h!]
	\begin{center}
		\tiny\begin{tabular}{|c|c|c|c|c|c|c|c|}
			\hline
			\multicolumn{1}{ |c| }{ $r  $ } &
			\multicolumn{1}{ |c| }{ G-ROM } &	
			\multicolumn{1}{ |c| }{ $\delta$ } &
			\multicolumn{1}{ |c| }{ DDC-ROM} &
			\multicolumn{1}{ |c| }{ $\delta $ } &
			\multicolumn{1}{ |c| }{ ICE-DDC-ROM } &			\multicolumn{1}{ |c| }{ $\delta$} &
			\multicolumn{1}{ |c| }{ CE-DDC-ROM } \\
			\hline
			2 & 1.78\textit{e}-03 & 1.00\textit{e}-05 & 1.32\textit{e}-03 & 1.00\textit{e}-04 & 1.46\textit{e}-05 & 1.00\textit{e}-04 & 1.48\textit{e}-05\\
			\hline
			3 & 1.93\textit{e}-04 & 1.00\textit{e}-05 & 1.55\textit{e}-04 & 
			1.00\textit{e}-06 &	 1.73\textit{e}-07 & 1.00\textit{e}-06 & 1.55\textit{e}-07\\
			\hline
			4 & 2.00\textit{e}-05 & 1.00\textit{e}-04 & 1.69\textit{e}-05 & 
			 1.00\textit{e}-06 & 9.46\textit{e}-08	& 1.00\textit{e}-06 & 9.49\textit{e}-08	\\
			\hline
			5 & 2.03\textit{e}-06 & 1.00\textit{e}-04 & 1.77\textit{e}-06 & 1.00\textit{e}-05 & 9.23\textit{e}-08 & 1.00\textit{e}-05 & 9.25\textit{e}-08 \\
			\hline
			6 & 2.23\textit{e}-07 & 1.00\textit{e}-05 & 2.02\textit{e}-07 & 1.00\textit{e}-05 & 9.37\textit{e}-08 & 1.00\textit{e}-05 & 9.38\textit{e}-08 \\
			\hline
			7 & 9.56\textit{e}-08 & 1.00\textit{e}-05 & 9.55\textit{e}-08 & 1.00\textit{e}-08 & 9.56\textit{e}-08 & 1.00\textit{e}-05 & 9.56\textit{e}-08 \\
			\hline		
		\end{tabular}
		\caption{Burgers equation, ROM differential filter, smooth initial condition, $\nu=10^{-1}$, and $d=7$: Average error in G-ROM, DDC-ROM, ICE-DDC-ROM, and CE-DDC-ROM for different $\delta$ and $r$ values.}\label{smooth-diff-effect}
	\end{center}
\end{table}

\begin{table}[h!]
	\begin{center}
		\smallskip
		\begin{tabular}{|c|c|c|c|c|}
			\hline
			\multicolumn{1}{ |c| }{ $r$ } &  
			\multicolumn{1}{ |c| }{ G-ROM } & 
			\multicolumn{1}{ |c| }{ DDC-ROM } & 
			\multicolumn{1}{ |c| }{ ICE-DDC-ROM }  &
			\multicolumn{1}{ |c| }{ CE-DDC-ROM }  \\
			\hline  
			2 &  1.78\textit{e}-03 & 1.32\textit{e}-03 & 1.46\textit{e}-05 & 1.48\textit{e}-05 \\
			\hline			
			3 &  1.93\textit{e}-04 & 1.55\textit{e}-04  & 1.73\textit{e}-07 & 1.55\textit{e}-07\\			
			\hline
			4 & 2.00\textit{e}-05 & 7.40\textit{e}-06 & 8.72\textit{e}-08  & 8.90\textit{e}-08	 \\ 
			\hline 
			5 &  2.03\textit{e}-06 & 7.29\textit{e}-07  & 9.23\textit{e}-08  & 9.25\textit{e}-08\\
			\hline 
			6 & 2.23\textit{e}-07 & 1.60\textit{e}-07   & 9.35\textit{e}-08 & 9.35\textit{e}-08   \\
			\hline
			7 & 9.56\textit{e}-08 & 9.56\textit{e}-08 & 9.56\textit{e}-08 & 9.56\textit{e}-08 \\
			\hline
		\end{tabular}
	\end{center}
	\caption{Burgers equation, ROM projection filter, smooth initial condition, $\nu=10^{-1}$, and $d=7$: Average error in G-ROM, DDC-ROM, ICE-DDC-ROM, and CE-DDC-ROM for different $r$ values. }\label{smooth-proj-effect} 
\end{table}

Next, we address (Q2), i.e., whether the CE has a significant effect on ROMs.
To this end, we test the DDC-ROM \eqref{ddc}, the ICE-DDC-ROM \eqref{ice-ddc-rom}, and the CE-DDC-ROM \eqref{ce-ddc}. 
We note that the ICE-DDC-ROM and CE-DDC-ROM include a representation of the Laplacian CE \eqref{comm_def2}, whereas the DDC-ROM does not. 
For various $r$ and $\delta$ values, we list the ROM error computed by using the ROM differential filter (Table \ref{smooth-diff-effect}) and the ROM projection filter  (Table \ref{smooth-proj-effect}).
We observe that the ICE-DDC-ROM and CE-DDC-ROM errors are consistently lower than the DDC-ROM error. 
We emphasize that, for low $r$ values, the ICE-DDC-ROM and CE-DDC-ROM errors are {\it two and even three orders of magnitude} lower than the DDC-ROM error.
Thus, we conclude that the CE plays a significant role in ROM development. 
Tables \ref{smooth-diff-effect} and \ref{smooth-proj-effect} also show that, as $r$ increases, the DDC-ROM error approaches the ICE-DDC-ROM and CE-DDC-ROM errors.

%


\subsubsection{\uppercase{Non-Smooth Initial Condition}}
\label{non-smooth-ic}

In this section, we investigate the effect of non-smooth initial conditions on the results obtained in Section~\ref{smooth-ic}.
To this end, we consider the following initial condition:
\begin{equation}
u_0(x)=\begin{cases}
\displaystyle~ 1, & x \in (0,1/2],\\
~\displaystyle 0, & x \in (1/2,1].
\end{cases}
\end{equation}

\begin{table}[h!]
	\begin{center}
		\begin{tabular}{|c|c|c|c|c|c|}
			\hline
			$r $  & $\delta $  & $||\mathcal{E}_{\Delta}[u_d]||_{L^2(L^2)}$  &			$r$  & $\delta $  & $||\mathcal{E}_{\Delta}[u_d]||_{L^2(L^2)}$   \\
			\hline
			2 &  1.00\textit{e}-01 & 1.85\textit{e}+00 & 9 &  1.00\textit{e}-01 & 2.13\textit{e}+00 \\
			\hline
			2 & 1.00\textit{e}-02 &1.12\textit{e}+00 & 9 & 1.00\textit{e}-02 & 2.25\textit{e}+01 \\
			\hline
			2 &  1.00\textit{e}-03 &1.90\textit{e}+00 &  9 &  1.00\textit{e}-03 & 7.12\textit{e}+02 \\
			\hline
			2 &  1.00\textit{e}-04 & 4.18\textit{e}+00 & 9 &  1.00\textit{e}-04 & 1.27\textit{e}+02\\
			\hline
			2 &  1.00\textit{e}-05 &4.29\textit{e}+00 & 9 &  1.00\textit{e}-05 & 1.17\textit{e}+02 \\
			\hline
			\hline
			3 &  1.00\textit{e}-01 &2.71\textit{e}+00 & 11 &  1.00\textit{e}-01 & 2.13\textit{e}+00 \\
			\hline
			3 & 1.00\textit{e}-02 & 6.70\textit{e}+00 & 11 & 1.00\textit{e}-02 & 2.14\textit{e}+01 \\
			\hline
			3 &  1.00\textit{e}-03 &1.07\textit{e}+01 & 11 &  1.00\textit{e}-03 & 7.14\textit{e}+02 \\
			\hline
			3 &  1.00\textit{e}-04 & 2.68\textit{e}+01 & 11 &  1.00\textit{e}-04 & 8.15\textit{e}+01 \\
			\hline
			3 &  1.00\textit{e}-05 & 2.75\textit{e}+01 &  11 &  1.00\textit{e}-05 & 3.59\textit{e}+01 \\
			\hline
			\hline                                   
			5 &  1.00\textit{e}-01 & 2.13\textit{e}+00 & 17 &  1.00\textit{e}-01 & 2.13\textit{e}+00 \\
			\hline
			5 & 1.00\textit{e}-02 & 5.00\textit{e}+01 & 17 & 1.00\textit{e}-02 & 2.13\textit{e}+01 \\
			\hline
			5 &  1.00\textit{e}-03 & 2.35\textit{e}+02 & 17 &  1.00\textit{e}-03 & 6.36\textit{e}+02 \\
			\hline
			5 &  1.00\textit{e}-04 & 3.15\textit{e}+02  & 17 &  1.00\textit{e}-04 & 1.10\textit{e}+02 \\
			\hline
			5 &  1.00\textit{e}-05 & 3.25\textit{e}+02 & 17 &  1.00\textit{e}-05 & 1.29\textit{e}+00 \\
			\hline
			\hline
			7 &  1.00\textit{e}-01 & 2.11\textit{e}+00 & 19 &  1.00\textit{e}-01 & 2.13\textit{e}+00 \\
			\hline
			7 & 1.00\textit{e}-02 & 2.90\textit{e}+01 & 19 & 1.00\textit{e}-02 & 2.13\textit{e}+01 \\
			\hline
			7 &  1.00\textit{e}-03 & 6.25\textit{e}+02 & 19 &  1.00\textit{e}-03 & 6.05\textit{e}+02  \\
			\hline
			7 &  1.00\textit{e}-04 & 2.95\textit{e}+02 & 19 &  1.00\textit{e}-04 & 1.20\textit{e}+02 \\
			\hline
			7 &  1.00\textit{e}-05 & 3.04\textit{e}+02 &  19 &  1.00\textit{e}-05 & 1.29\textit{e}+00 \\
			\hline
		\end{tabular}
	\end{center}
	\caption{Burgers equation, ROM differential filter, non-smooth initial condition, $\nu=10^{-1}$, and $d=19$: Average CE for different $\delta$ and $r$ values.}\label{step-diff}
\end{table}

\begin{table}[h!]
	\begin{center} 
		\tiny\begin{tabular}{|l|c|c|c|c|c|c|r|}
			\hline
			$r$ &  2 & 3 & 5 & 7 & 9 & 11 & 17  \\ \hline
			$ \bold ||\mathcal{E}_{\Delta}[u_d]||_{L^2(L^2)}$ & 4.29\textit{e}+00 & 2.75\textit{e}+01 & 3.25\textit{e}+02 & 3.04\textit{e}+02 & 1.17\textit{e}+02 & 3.60\textit{e}+01 &  5.27\textit{e}-01  \\ \hline		
		\end{tabular}
	\end{center}
	\caption{Burgers equation, ROM projection filter, non-smooth initial condition, $\nu=10^{-1}$, and $d=19$: Average CE for different $r$ values.}\label{step-proj}
\end{table}

As in Section~\ref{smooth-ic}, we start by addressing (Q1), i.e., whether the CE exists.
For various $r$ and $\delta$ values, we list the CE computed by using the ROM differential filter (Table \ref{step-diff}) and the ROM projection filter (Table \ref{step-proj}).
We draw the same main conclusion as in Section~\ref{smooth-ic}: 
The CE exists for both filters.
This time, however, as $r$ increases, the CE does not decrease.
Furthermore, for the ROM differential filter (Table \ref{step-diff}), we note that for low $r$ values, as $\delta$ decreases, the CE increases.

\begin{table}[h!]
	\begin{center}
		\tiny\begin{tabular}{|c|c|c|c|c|c|c|c|}
			\hline
			\multicolumn{1}{ |c| }{ $r$ } & 
			\multicolumn{1}{ |c| }{G-ROM } & 
			\multicolumn{1}{ |c| }{ $\delta$ } &
			\multicolumn{1}{ |c| }{DDC-ROM } & 
			\multicolumn{1}{ |c| }{  $\delta$ } &
			\multicolumn{1}{ |c| }{ ICE-DDC-ROM } &
			\multicolumn{1}{ |c| }{ $\delta$ } &  
			\multicolumn{1}{ |c| }{ CE-DDC-ROM }  \\
			\hline 
			2 & 7.27\textit{e}-03 & 1.00\textit{e}-05 & 6.93\textit{e}-03  & 1.00\textit{e}-03 & 2.67\textit{e}-04  & 1.00\textit{e}-02 & 1.07\textit{e}-03\\
			\hline
			3 & 1.51\textit{e}-02 & 1.00\textit{e}-04 & 8.16\textit{e}-03 & 1.00\textit{e}-04 & 5.85\textit{e}-05	& 1.00\textit{e}-02 & 1.36\textit{e}-03	\\
			\hline
			5 & 4.22\textit{e}-03 & 1.00\textit{e}-04 & 4.22\textit{e}-03 & 1.00\textit{e}-07 & 5.97\textit{e}-07 & 1.00\textit{e}-04 & 3.13\textit{e}-04 \\
			\hline
			7 & 9.59\textit{e}-04 & 1.00\textit{e}-04 & 9.59\textit{e}-04 & 1.00\textit{e}-07 & 1.32\textit{e}-07  & 1.00\textit{e}-04 & 6.27\textit{e}-05 \\
			\hline
			9 & 	2.34\textit{e}-04 & 1.00\textit{e}-04 & 2.35\textit{e}-04 & 1.00\textit{e}-06 & 1.42\textit{e}-07  &  1.00\textit{e}-05  & 2.07\textit{e}-06 \\
			\hline
			11 & 5.44\textit{e}-05 & 1.00\textit{e}-07 & 5.56\textit{e}-05 & 1.00\textit{e}-07 & 1.38\textit{e}-07 & 1.00\textit{e}-06 & 1.48\textit{e}-07 \\
			\hline
			17 & 4.97\textit{e}-07 & 1.00\textit{e}-07 & 5.06\textit{e}-07 & 1.00\textit{e}-07 & 1.35\textit{e}-07 & 1.00\textit{e}-07 & 1.19\textit{e}-07 \\
			\hline
			19 & 1.30\textit{e}-07 & 1.00\textit{e}-07 & 1.30\textit{e}-07 & 1.00\textit{e}-07 & 1.30\textit{e}-07 & 1.00\textit{e}-07 & 1.30\textit{e}-07  \\
			\hline          			
		\end{tabular}
	\end{center}
	\caption{Burgers equation, ROM differential filter, non-smooth initial condition, $\nu=10^{-1}$, and $d=19$: Average error in G-ROM, DDC-ROM, ICE-DDC-ROM, and CE-DDC-ROM for different $\delta$ and $r$ values.}	
	\label{step-diff-effect}
\end{table}	

\begin{table}[h!]
	\begin{center}
		\smallskip
		\begin{tabular}{|c|c|c|c|c|}
			\hline
			\multicolumn{1}{ |c| }{ $r$ } & 
			\multicolumn{1}{ |c| }{ G-ROM } & 
			\multicolumn{1}{ |c| }{ DDC-ROM } & 
			\multicolumn{1}{ |c| }{ ICE-DDC-ROM } & 
			\multicolumn{1}{ |c| }{ CE-DDC-ROM }  \\
			\hline 
			2 & 7.27\textit{e}-03 & 6.93\textit{e}-03 & 2.71\textit{e}-04  & 3.71e-03 \\
			\hline
			3 & 1.51\textit{e}-02 & 8.14\textit{e}-03 & 5.89\textit{e}-05 & 3.84\textit{e}-03 \\
			\hline
			5 & 4.22\textit{e}-03 & 4.22\textit{e}-03 & 5.97\textit{e}-07 & 3.14\textit{e}-04 \\
			\hline
			7 & 9.59\textit{e}-04 & 9.59\textit{e}-04 & 1.32\textit{e}-07 & 6.44\textit{e}-05 \\
			\hline
			9 & 2.34\textit{e}-04 & 2.34\textit{e}-04 & 1.39\textit{e}-07 & 4.26\textit{e}-06 \\
			\hline		
			11 & 5.44\textit{e}-05 & 5.44\textit{e}-05 & 1.30\textit{e}-07 & 1.28\textit{e}-07 \\
			\hline
			17 & 4.97\textit{e}-07 & 5.00\textit{e}-07 & 1.34\textit{e}-07 & 1.15\textit{e}-07 \\
			\hline
			19 & 1.30\textit{e}-07 & 1.30\textit{e}-07 & 1.30\textit{e}-07 & 1.30\textit{e}-07 \\
			\hline          			
		\end{tabular}
	\end{center}
	\caption{Burgers equation, ROM projection filter, non-smooth initial condition, $\nu=10^{-1}$, and $d=19$: Average error in G-ROM, DDC-ROM, ICE-DDC-ROM, and CE-DDC-ROM for different $r$ values.}	
	\label{step-proj-effect}	
\end{table}	

Next, we address (Q2), i.e., whether the CE has a significant effect on ROMs.
As in Section~\ref{smooth-ic}, we test the DDC-ROM \eqref{ddc}, the ICE-DDC-ROM \eqref{ice-ddc-rom}, and the CE-DDC-ROM \eqref{ce-ddc}. 
We note again that the ICE-DDC-ROM and CE-DDC-ROM include a representation of the Laplacian CE \eqref{comm_def2}, whereas the DDC-ROM does not. 
For various $r$ and $\delta$ values, we list the ROM error computed by using the ROM differential filter (Table \ref{step-diff-effect}) and the ROM projection filter  (Table \ref{step-proj-effect}).
We draw the same main conclusion as in Section~\ref{smooth-ic}: 
The ICE-DDC-ROM and CE-DDC-ROM errors are consistently lower than the DDC-ROM error. 
Furthermore, for low $r$ values, the ICE-DDC-ROM and CE-DDC-ROM errors are {\it one and even two orders of magnitude} lower than the DDC-ROM error.
Thus, we conclude again that the CE plays a significant role in ROM development. 
Tables \ref{step-diff-effect} and \ref{step-proj-effect} also show that, as $r$ increases, the DDC-ROM error approaches the ICE-DDC-ROM and CE-DDC-ROM errors.

\subsubsection{\uppercase{Lower Viscosity} ($\nu=10^{-3}$)}
	\label{sec:burgers-lower-viscosity}

In this section, we investigate the effect of lower viscosity on the results obtained in Sections~\ref{smooth-ic} and~\ref{non-smooth-ic}.
For clarity, we only present results for the smooth initial conditions used in Section~\ref{smooth-ic}; the results for the non-smooth initial condition used in Section~\ref{non-smooth-ic} were similar.
As in the previous sections, we start by addressing (Q1), i.e., whether the CE exists.
For various $r$ and $\delta$ values, we list the CE computed by using the ROM differential filter (Table \ref{smooth_ce-re-1000-diff-burger}) and the ROM projection filter (Table \ref{tab:burgers-ce-proj-low-nu}).
We draw the same main conclusion as in Sections~\ref{smooth-ic} and~\ref{non-smooth-ic}: 
The CE exists for both filters.

\begin{table}[h]
	\begin{center}
		\begin{tabular}{|c|c|c|}
			\hline
			\multicolumn{1}{ |c| }{ $r$ } &
			\multicolumn{1}{ |c| }{ $   \delta $ } &		
			\multicolumn{1}{ |c| }{ $  ||\mathcal{E}_{\Delta}[u_d]||_{L^2(L^2)}  $ } \\
			\hline
			1 &  1.00\textit{e}-01 & 3.61\textit{e}-04  \\ 
			\hline
			1 &  1.00\textit{e}-02 & 4.39\textit{e}-04  \\ 
			\hline
			1 & 1.00\textit{e}-04 & 4.40\textit{e}-04  \\ 
			\hline
			1 &  1.00\textit{e}-06 & 4.40\textit{e}-04  \\ 
			\hline
			\hline
			2 & 1.00\textit{e}-01 & 2.09\textit{e}-06  \\
			\hline 
			2 & 1.00\textit{e}-02 & 4.28\textit{e}-06  \\
			\hline
			2 & 1.00\textit{e}-04 & 4.32\textit{e}-06  \\
			\hline
			2 & 1.00\textit{e}-06 & 4.32\textit{e}-06  \\
			\hline
			\hline						
			3 &  1.00\textit{e}-01 & 1.36\textit{e}-15  \\
			\hline
			3 &  1.00\textit{e}-02& 4.80\textit{e}-17  \\
			\hline
			3 &  1.00\textit{e}-04 & 3.80\textit{e}-17  \\
			\hline
			3 &  1.00\textit{e}-06 & 3.74\textit{e}-17  \\
			\hline            
		\end{tabular}
		\caption{Burgers equation, ROM differential filter, smooth initial condition, $\nu =10^{-3}$, and $d=3$: Average CE for different $\delta$ and $r$ values.
		\label{smooth_ce-re-1000-diff-burger}
		}
	\end{center}
\end{table}

\begin{table}[h]
	\begin{center}
		\begin{tabular}{|c|c|c|c|}
			\hline
			$r $ &  1 & 2 & 3 \\
			\hline
			$||\mathcal{E}_{\Delta}[u_d]||_{L^2(L^2)}$    & 4.40\textit{e}-04  & 4.32\textit{e}-06  & 0 \\
			\hline
		\end{tabular}
		\caption{Burgers equation, ROM projection filter, smooth initial condition, $\nu=10^{-3}$, and $d=3$: Average CE for different $r$ values.
		\label{tab:burgers-ce-proj-low-nu}
		}
	\end{center}
\end{table}

Next, we address (Q2), i.e., whether the CE has a significant effect on ROMs.
As in Sections~\ref{smooth-ic} and~\ref{non-smooth-ic}, we test the DDC-ROM \eqref{ddc}, the ICE-DDC-ROM \eqref{ice-ddc-rom}, and the CE-DDC-ROM \eqref{ce-ddc}. 
We note again that the ICE-DDC-ROM and CE-DDC-ROM include a representation of the Laplacian CE \eqref{comm_def2}, whereas the DDC-ROM does not. 
For various $r$ and $\delta$ values, we list the ROM error computed by using the ROM differential filter (Table \ref{burgers-smooth-diff-effect-re-1000}) and the ROM projection filter  (Table \ref{soom-re-1000-proj-burger}).
We draw the same main conclusion as in Sections~\ref{smooth-ic} and~\ref{non-smooth-ic}: 
The ICE-DDC-ROM and CE-DDC-ROM errors are consistently lower than the DDC-ROM error. 
Furthermore, for low $r$ values, the ICE-DDC-ROM and CE-DDC-ROM errors are {\it two and even three orders of magnitude} lower than the DDC-ROM error.
Thus, we conclude again that the CE plays a significant role in ROM development. 
We also note that, for the Burgers equation, the lower viscosity ($\nu=10^{-3}$) results are similar to the higher viscosity ($\nu=10^{-1}$) results. 

 \begin{table}[h]
	\begin{center}
		\tiny\begin{tabular}{|c|c|c|c|c|c|c|c|}
			\hline
			\multicolumn{1}{ |c| }{ $ r $ } & 
			\multicolumn{1}{ |c| }{ G-ROM } & 
			\multicolumn{1}{ |c| }{ $ \delta $ } &
			\multicolumn{1}{ |c| }{  DDC-ROM  } & 
			\multicolumn{1}{ |c| }{ $  \delta$ } &
			\multicolumn{1}{ |c| }{  ICE-DDC-ROM  } & 
			\multicolumn{1}{ |c| }{ $  \delta $ } &
			\multicolumn{1}{ |c| }{ CE-DDC-ROM }  \\
			\hline 
			1 & 2.21\textit{e}-07 & 1.00\textit{e}-03 & 5.81\textit{e}-08  & 1.00\textit{e}-05 & 3.50\textit{e}-10  & 1.00\textit{e}-06 &  4.07\textit{e}-10 \\
			\hline
			2 & 9.33\textit{e}-10  & 1.00\textit{e}-08  & 6.68\textit{e}-10 & 1.00\textit{e}-08 & 5.35\textit{e}-13 &1.00\textit{e}-08 & 5.37\textit{e}-13
			\\          \hline  
			3 & 3.46\textit{e}-12 & 1.00\textit{e}-08  & 3.46\textit{e}-12 & 1.00\textit{e}-08 & 3.46\textit{e}-12 & 1.00\textit{e}-08 & 3.46\textit{e}-12 
			\\          \hline 
		\end{tabular}
	\end{center}		
    \caption{Burgers equation, ROM differential filter, smooth initial condition, $\nu=10^{-3}$, and $d=3$: Average error in G-ROM, DDC-ROM, ICE-DDC-ROM, and CE-DDC-ROM for different $\delta$ and $r$  values.
    \label{burgers-smooth-diff-effect-re-1000}
    }
\end{table}

\begin{table}[h]
	\begin{center}
		\smallskip
		\begin{tabular}{|c|c|c|c|c|}
			\hline
			\multicolumn{1}{ |c| }{ $r$ } & 
			\multicolumn{1}{ |c| }{ G-ROM} & 
			\multicolumn{1}{ |c| }{ DDC-ROM } & 
			\multicolumn{1}{ |c| }{ ICE-DDC-ROM } & 
			\multicolumn{1}{ |c| }{ CE-DDC-ROM }  \\
			\hline 
			1 & 2.21\textit{e}-07 & 1.35\textit{e}-07 & 3.49\textit{e}-10  & 4.07\textit{e}-10 \\
			\hline
			2 & 9.33\textit{e}-10  & 6.68\textit{e}-10  & 5.35\textit{e}-13 &  5.37\textit{e}-13
			\\          \hline  
			3 & 3.46\textit{e}-12 & 3.46\textit{e}-12  & 3.46\textit{e}-12  & 3.46\textit{e}-12 
			\\          \hline 
		\end{tabular}
	\end{center}		
	\caption{Burgers equation, ROM projection filter, smooth initial condition, $\nu=10^{-3}$, and $d=3$:  Average error in G-ROM, DDC-ROM, ICE-DDC-ROM, and CE-DDC-ROM for different $r$ values.\label{soom-re-1000-proj-burger}
	}
\end{table}

\subsection{\uppercase{Experiment 2: Flow past a circular cylinder}}
\label{sec:test-problem-setup}

\begin{figure}[h!]
	\begin{center}
		\includegraphics[width=0.7\textwidth,height=0.24\textwidth, trim=0 0 0 0, clip]{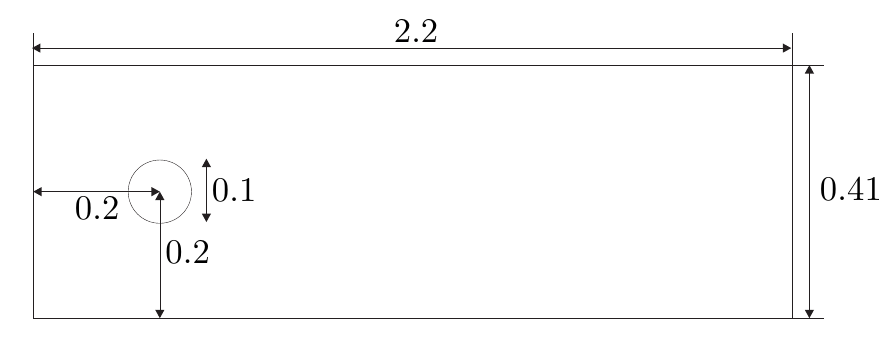}
	\end{center}
	\caption{\label{cyldomain} Channel flow around a cylinder domain.}
\end{figure}

In our second experiment, we consider a 2D channel flow past a circular cylinder. The domain is a $2.2\times 0.41$ rectangular channel with a radius $=0.05$ cylinder, centered at $(0.2,0.2)$, see Fig.~\ref{cyldomain}.  
No slip boundary conditions are prescribed on the walls and cylinder, and the inflow and outflow profiles are given by~\cite{mohebujjaman2018physically,mohebujjaman2017energy}
\begin{align*}
u_{1}(0,y,t)&=u_{1}(2.2,y,t)=\frac{6}{0.41^{2}}y(0.41-y), \\ u_{2}(0,y,t)&=u_{2}(2.2,y,t)=0,
\end{align*} 
where $\bu=\langle u_1, u_2 \rangle$.  Even though the parabolic outflow condition is not physical, the Dirichlet boundary condition is often used for ROMs, since it reduces the theoretical and computational complexities. There is no forcing $(f=0)$ and the flow starts from rest. We run the DNS of the NSE~\eqref{eqn:nse-1}-\eqref{eqn:nse-2} from rest ($t=0$) until the simulation time $t=17$.  
We use the point-wise divergence-free, LBB stable $(P_2, P_1^{disc})$ Scott-Vogelius FE pair on a barycenter refined regular triangle mesh. The mesh provides $16178$ velocity, and $11907$ pressure degrees of freedom. The time step size $\Delta t=0.002$ is used for both DNS and ROM time evolution. 
We utilize the commonly used linearized BDF2 temporal discretization, together with the FE spatial discretization. On the first time step, we use a backward Euler scheme so that we have two initial time step solutions required for the BDF2 scheme.
The scheme for $n=1,2,\cdots,$ is:  Find $(\bu_h^{n+1}, p_h^{n+1})\in (X^h, Q^h)$ satisfying for every $(\bv_h, q_h)\in (X^h, Q^h)$,
\begin{eqnarray}
\left(\frac{3\bu_h^{n+1}-4\bu_h^n+\bu_h^{n-1}}{2\Delta t}, \bv_h\right)+((2\bu_h^{n}-\bu_h^{n-1})\cdot\nabla \bu_h^{n+1}, \bv_h)\nonumber\\-(p_h^{n+1}, \nabla\cdot \bv_h)+\nu(\nabla \bu_h^{n+1},\nabla \bv_h) = 0, \label{disc-nse}\\
(\nabla\cdot u_h, q_h)=0.
\end{eqnarray}
We collect $2500$ snapshots at each time step from $t=5$ to $t=10$ to construct the ROM basis.

\subsubsection{\uppercase{Reynolds Number} $Re=1$}
	\label{sec:re-1}

In this section, we present numerical results for $Re=1$, which corresponds to $\nu=0.1$. 
As in Section~\ref{burgers}, we address questions (Q1) and (Q2).
We start with (Q1), i.e., whether the CE exists.
For various $r$ and $\delta$ values, we list the CE computed by using the ROM differential filter (Table \ref{ceRe_1_diff}) and the ROM projection filter  (Table \ref{ceproj_re_1_1}).
The main conclusion is that the CE exists for both filters, just as for the Burgers equation in Section~\ref{burgers}.
However, the CE in Tables \ref{ceRe_1_diff} and \ref{ceproj_re_1_1} is much lower than the CE for Burgers equation in Section~\ref{burgers}.

\begin{table}[h!]
	\begin{center} 
		\scriptsize\begin{tabular}{|l|c|c|c|c|r|}
			\hline
			$r$& $\delta=$1.00\textit{e}-04 & $\delta=$1.00\textit{e}-03 & $\delta=$1.00\textit{e}-2 &$\delta=$1.00\textit{e}-01 &$\delta=$1.50\textit{e}-01  \\ \hline
			3 & 4.84\textit{e}-11 & 2.53\textit{e}-11 & 1.01\textit{e}-12 & 1.51\textit{e}-14 & 8.81\textit{e}-15\\ \hline
			4 & 4.83\textit{e}-11 & 2.53\textit{e}-11 & 1.01\textit{e}-12 & 2.12\textit{e}-14 & 8.46\textit{e}-15\\ \hline
			7 & 5.27\textit{e}-11 & 1.91\textit{e}-11 & 6.74\textit{e}-13 & 3.02\textit{e}-13 & 1.88\textit{e}-13\\ \hline
			9 & 3.67\textit{e}-11 & 6.04\textit{e}-12 & 6.98\textit{e}-13 & 1.04\textit{e}-12 & 6.83\textit{e}-13\\ \hline
			11& 3.72\textit{e}-11 & 6.15\textit{e}-12 & 6.94\textit{e}-13 & 1.05\textit{e}-12 & 6.05\textit{e}-13\\ \hline
			13& 3.68\textit{e}-11 & 5.97\textit{e}-12 & 6.95\textit{e}-13 & 1.01e-12 & 6.10e-13\\ \hline
			15& 3.71\textit{e}-11 & 6.90\textit{e}-12 & 6.68\textit{e}-13 & 1.12\textit{e}-12 & 7.23\textit{e}-13\\ \hline
		\end{tabular}
	\end{center}
	\caption{NSE, $Re=1$, ROM differential filter, $d=16$: Average CE for different $\delta$ and $r$ values.}\label{ceRe_1_diff}
\end{table}

\begin{table}[h!]
	\begin{center} 
		\begin{tabular}{|l|c|c|c|c|c|r|}
			\hline
			$r$ &  2 & 3 & 4 & 5 & 6 & 7\\ \hline
			$\bold ||\mathcal{E}_{\Delta}[\bu_d]||_{L^2(L^2)}$  & 6.30\textit{e}-13 & 7.41\textit{e}-11 & 7.45\textit{e}-11 & 7.31\textit{e}-11 & 4.00\textit{e}-11 & 2.84\textit{e}-11 \\ \hline
		\end{tabular}
	\end{center}
	\caption{NSE, $Re=1$, ROM projection filter, and $d=7$: Average CE for different $r$ values.}\label{ceproj_re_1_1}
\end{table}

\begin{table}[h!]
	\begin{center} 
		\begin{tabular}{|l|c|c|c|r|}
			\hline
			$r$  & G-ROM & DDC-ROM & ICE-DDC-ROM & CE-DDC-ROM  \\ \hline
			3 & 5.51\textit{e}-05& 5.51\textit{e}-05 &5.51\textit{e}-05 & 5.12\textit{e}-05 \\ \hline
			4 &  5.50\textit{e}-03 & 5.50\textit{e}-03& 5.50\textit{e}-03 & 1.14\textit{e}-04 \\ \hline	
			5 & 1.00\textit{e}-02 & 1.00\textit{e}-02& 1.00\textit{e}-02 & 6.90\textit{e}-03 \\ \hline	
		\end{tabular}
	\end{center}
	\caption{NSE, $Re=1$, $d=7$, ROM projection filter: Average error on G-ROM, DDC-ROM, ICE-DDC-ROM, and CE-DDC-ROM for different $r$ values. 
	}\label{proj_re_1_effect}
\end{table}

\begin{table}[h!]
	\begin{center} 
		\begin{tabular}{|l|c|c|c|r|}
			\hline
			$r$&$\delta$   & DDC-ROM & ICE-DDC-ROM & CE-DDC-ROM  \\ \hline
			3 & 1.00\textit{e}-04& 5.51\textit{e}-05 &5.51\textit{e}-05 & 5.12\textit{e}-05 \\ \hline
			3 & 1.00\textit{e}-03& 5.51\textit{e}-05 & 5.51\textit{e}-05& 5.23\textit{e}-05 \\ \hline
			3 & 1.00\textit{e}-02& 5.51\textit{e}-05 & 5.51\textit{e}-05& 5.51\textit{e}-05 \\ \hline		
			4 & 1.00\textit{e}-04 &5.50\textit{e}-03&5.50\textit{e}-03& 1.14\textit{e}-04\\ \hline
			4 & 1.00\textit{e}-03 &5.50\textit{e}-03&5.50\textit{e}-03& 1.79\textit{e}-04\\ \hline
			4 & 1.00\textit{e}-02 &5.50\textit{e}-03&5.50\textit{e}-03& 3.30\textit{e}-03\\ \hline
			5 & 1.00\textit{e}-04&1.00\textit{e}-02 & 1.00\textit{e}-02& 5.90\textit{e}-03\\ \hline
			5 & 1.00\textit{e}-03&1.00\textit{e}-02 & 1.00\textit{e}-02& 7.00\textit{e}-03\\ \hline
			5 & 1.00\textit{e}-02&1.00\textit{e}-02 & 1.00\textit{e}-02& 9.90\textit{e}-03\\ \hline	
		\end{tabular}
	\end{center}
	\caption{NSE, $Re=1$, $d=7$, ROM differential filter: Average error in G-ROM, DDC-ROM, ICE-DDC-ROM, and CE-DDC-ROM for different $\delta$ and $r$ values. 
	\label{diff_re_1_ROM_effect}
	}
\end{table}

Next, we address (Q2), i.e., whether the CE has a significant effect on ROMs.
To this end, we test the DDC-ROM \eqref{ddc}, the ICE-DDC-ROM \eqref{ice-ddc-rom}, and the CE-DDC-ROM \eqref{ce-ddc}. 
We note that the ICE-DDC-ROM and CE-DDC-ROM include a representation of the Laplacian CE \eqref{comm_def2}, whereas the DDC-ROM does not. 
For various $r$ and $\delta$ values, we list the ROM error computed by using the ROM differential filter (Table \ref{diff_re_1_ROM_effect}) and the ROM projection filter  (Table \ref{proj_re_1_effect}).
We observe that the ICE-DDC-ROM and CE-DDC-ROM errors are lower than the DDC-ROM error. 
Thus, we conclude that the CE plays a significant role in ROM development. 
We note, however, that the ICE-DDC-ROM and CE-DDC-ROM errors are only a factor of $2$ lower than the DDC-ROM error.
This is different from the Burgers equation in Section~\ref{burgers}, where the ICE-DDC-ROM and CE-DDC-ROM errors were orders of magnitude lower than the DDC-ROM error.
We also note that the DDC-ROM and ICE-DDC-ROM errors in Tables~\ref{proj_re_1_effect} and \ref{diff_re_1_ROM_effect} are the same.  
This is due to the fact that the magnitude of the CE is of $\mathcal{O}(10^{-9})$, which is much lower than the magnitude of the other ROM terms (i.e., $\mathcal{O}(10^{-1})$).
On the other hand, finding a data-driven model for the CE yields terms of the same order of magnitude as  the ROM terms, i.e., $\mathcal{O}(10^{-1})$.
This explains why the CE-DDC-ROM error is significantly smaller than, say, the ICE-DDC-ROM error.

\subsubsection{\uppercase{Reynolds Number} $Re=100$}\label{re-100:case}
To investigate the effect of the Reynolds number on the results in Section~\ref{sec:re-1}, in this section we consider $Re=100$, which corresponds to $\nu=10^{-3}$. 
To generate the ROM basis, we collect $166$ snapshots, which are the FE solutions at each time step from $t=7$ to $t=7.332$. 

As in Section~\ref{sec:re-1}, we start by addressing (Q1), i.e., whether the CE exists.
For various $r$ and $\delta$ values, we list the CE computed by using the ROM differential filter (Table \ref{cediff}) and the ROM projection filter (Table \ref{ceproj1}).
As in Section~\ref{sec:re-1}, we observe that the CE exists for both filters.
We also note that the magnitude of the CE for $Re=100$ is much higher than the magnitude of the CE for $Re=1$.

\begin{table}[h!]
	\begin{center} 
		\begin{tabular}{|l|c|c|c|c|r|}
			\hline
			$r$& $\delta=$1.00\textit{e}-04 & $\delta=$1.00\textit{e}-03 & $\delta=$1.00\textit{e}-2 &$\delta=$1.00\textit{e}-01 &$\delta=$1.50\textit{e}-01  \\ \hline
			3 & 9.84\textit{e}-01& 9.83\textit{e}-01 & 9.22\textit{e}-01& 1.68\textit{e}-01& 8.81\textit{e}-02 \\ \hline
			4 & 1.86\textit{e}+00& 1.85\textit{e}+00 & 1.54\textit{e}+00 & 1.13\textit{e}-01 & 5.32\textit{e}-02 \\ \hline
			7 & 3.06\textit{e}+00& 3.05\textit{e}+00 & 2.35\textit{e}+00& 1.01\textit{e}-01& 4.59\textit{e}-02\\ \hline
			9 & 3.94\textit{e}+00& 3.92\textit{e}+00 & 2.65\textit{e}+00& 9.68\textit{e}-02& 4.38\textit{e}-02\\ \hline
			11& 2.15\textit{e}+00& 2.13\textit{e}+00 & 1.19\textit{e}+00& 2.93\textit{e}-02& 1.32\textit{e}-02 \\ \hline
			13& 1.90\textit{e}+00& 1.88\textit{e}+00 & 9.75\textit{e}-01& 2.31\textit{e}-02 & 1.04\textit{e}-02 \\ \hline
			15& 7.80\textit{e}-01 & 7.68\textit{e}-01 &3.27\textit{e}-01 & 6.20\textit{e}-03& 2.80\textit{e}-03 \\ \hline
		\end{tabular}
	\end{center}
	\caption{NSE, $Re=100$, ROM differential filter, $d=16$: Average CE for different $\delta$ and $r$ values.}\label{cediff}
\end{table}

\begin{table}[h!]
	\begin{center} 
		\begin{tabular}{|l|c|c|c|c|c|r|}
			\hline
			$r$ &  2 & 3 & 4 & 5 & 6 & 7\\ \hline
			\tiny{$\bold ||\mathcal{E}_{\Delta}[\bu_d]||_{L^2(L^2)}$}  & 2.05\textit{e}+00 & 9.70\textit{e}-01 & 1.75\textit{e}+00 & 2.26\textit{e}+00 & 1.65\textit{e}+00 & 4.46\textit{e}-14 \\ \hline
		\end{tabular}
	\end{center}
	\caption{NSE, $Re=100$, ROM projection filter, and $d=7$: Average CE of different $r$ values.}\label{ceproj1}
\end{table}

Next, we address (Q2), i.e., whether the CE has a significant effect on ROMs.
To this end, we test the DDC-ROM \eqref{ddc}, the ICE-DDC-ROM \eqref{ice-ddc-rom}, and the CE-DDC-ROM \eqref{ce-ddc}. 
We note that the ICE-DDC-ROM and CE-DDC-ROM include a representation of the Laplacian CE \eqref{comm_def2}, whereas the DDC-ROM does not. 
For various $r$ and $\delta$ values, we list the ROM error computed by using the ROM differential filter (Table \ref{diff_re_100_effect}) and the ROM projection filter  (Table \ref{proj_re_100_effect}).
We observe that the ICE-DDC-ROM and CE-DDC-ROM errors are lower than the DDC-ROM error, just as in the previous cases. 
This time, however, the improvements in the ICE-DDC-ROM and CE-DDC-ROM over the DDC-ROM are small. 
Thus, for the $Re=100$ case, we conclude that the CE plays only a minor role in ROM development. 
This happens because, in this case, the magnitude of the CE is much lower than the magnitude of the Correction term.
The explanation for this difference is that, in~\eqref{ddc-ce}, the CE is multiplied by $\nu$, the viscosity coefficient.
Thus, for higher $Re$, the CE (which arises from the diffusion term) is dominated by the Correction term (which arises from the nonlinear term).
This explains why the CE is important for low $Re$ flows, but not for higher $Re$ flows.

\begin{table}[h!]
	\begin{center} 
		\begin{tabular}{|l|c|c|c|c|r|}
			\hline
			$r$& $d$  & G-ROM & DDC-ROM & ICE-DDC-ROM & CE-DDC-ROM  \\ \hline
			3 & 7   & 8.70\textit{e}-03 & 3.92\textit{e}-04 & 3.68\textit{e}-04  &   3.70\textit{e}-04  \\ \hline
			4 & 9   & 3.57\textit{e}-02 &  3.61\textit{e}-04 &  3.52\textit{e}-04 & 3.55\textit{e}-04     \\ \hline
			5 & 9   & 2.17\textit{e}-02 &  4.06\textit{e}-04 & 3.25\textit{e}-04 &   3.28\textit{e}-04   \\ \hline
		\end{tabular}
	\end{center}
	\caption{NSE, $Re=100$, ROM projection filter: Average error in G-ROM, DDC-ROM, ICE-DDC-ROM, and CE-DDC-ROM for different $d$ and $r$ values. 
	}\label{proj_re_100_effect}
\end{table}

\begin{table}[h!]
	\begin{center} 
		\begin{tabular}{|l|c|c|c|c|r|}
			\hline
			$r$& $d$& $\delta$        & DDC-ROM          & ICE-DDC-ROM        & CE-DDC-ROM  \\ \hline
			3 & 7   & 1.00\textit{e}-3& 3.92\textit{e}-04& 3.68\textit{e}-04  & 3.70\textit{e}-04  \\ \hline
			3 & 7   & 1.00\textit{e}-2& 3.92\textit{e}-04& 3.70\textit{e}-04  & 3.71\textit{e}-04  \\ \hline
			4 & 9   & 1.00\textit{e}-3& 3.61\textit{e}-04& 3.52\textit{e}-04  & 3.55\textit{e}-04 \\ \hline
			4 & 9   & 1.00\textit{e}-2& 3.61\textit{e}-04& 3.51\textit{e}-04  & 3.53\textit{e}-04   \\ \hline
			5 & 9   & 1.00\textit{e}-3& 4.06\textit{e}-04& 3.23\textit{e}-04  & 3.82\textit{e}-04  \\ \hline
			5 & 9   & 1.00\textit{e}-3& 4.06\textit{e}-04& 3.32\textit{e}-04  & 3.80\textit{e}-04  \\ \hline
		\end{tabular}
	\end{center}
	\caption{NSE, $Re=100$, ROM differential filter: Average error in G-ROM, DDC-ROM, ICE-DDC-ROM, and CE-DDC-ROM for different $\delta$ and $r$ values. 
	}\label{diff_re_100_effect}
\end{table}

\section{\uppercase{Conclusions and Future Work}}\label{conclusion}

In this paper, we investigated theoretically and computationally whether the {\it commutation error (CE)} exists, i.e., whether differentiation and ROM spatial filtering commute.
To our knowledge, this is the first investigation of the CE in a ROM context.
We studied whether there is a CE for the Laplacian for two ROM filters: the ROM projection and the ROM differential filter.
Furthermore, when the CE was nonzero, we investigated whether it had any significant effect on the ROM development.
To this end, we considered the {\it data-driven correction ROM (DDC-ROM)}~\cite{xie2018data}, in which the Correction term (which is generally added to improve the ROM's accuracy) is modeled by using the available data.
To investigate the effect of the CE on the DDC-ROM, we considered the commutation error DDC-ROM (CE-DDC-ROM), in which available data is used to model not only the Correction term, but also the CE.
Finally, we also used the ideal CE-DDC-ROM (ICE-DDC-ROM), which is the DDC-ROM supplemented with a fine resolution representation of the CE. 
When the CE-DDC-ROM and ICE-DDC-ROM yielded more accurate results than the standard DDC-ROM, we concluded that the CE has a significant effect on the ROM development.
As numerical tests, we used the Burgers equation with viscosities $\nu=10^{-1}$ and $\nu=10^{-3}$ and a 2D flow past a circular cylinder at Reynolds numbers $Re=1$ and $Re=100$.
For the Burgers equation test case, we considered smooth and non-smooth initial conditions.

The most important conclusions of our theoretical and numerical investigation are the following: 
(i) The CE exists for all cases considered.
(ii) The CE has a significant effect on the ROM development for low Reynolds numbers, but not so much for higher Reynolds numbers.
This happens because, for higher Reynolds numbers, the CE (which arises from the diffusion term) is dominated by the Correction term (which arises from the nonlinear term).
We note that, for the Burgers equation, the CE had a significant effect on ROMs even for low viscosity values; however, for non-smooth initial conditions (results not included), the CE effect was lower than the CE effect for higher viscosity.
(iii) The non-smooth initial conditions (in the Burgers equation) decreased the effect of the CE on the ROM development.

These first steps in the theoretical and numerical investigation of the CE showed that, in some cases, it can be significant and has to be modeled.
There are, however, several other research directions that need to be pursued for a better understanding of the ROM CE.
For example, we plan to investigate whether there is an upper bound for the Reynolds number for which the CE has a significant effect on the ROM. 
Furthermore, we plan to study the ROM CE for differential operators that are different from the Laplacian, e.g., first-order spatial derivatives, such as those in the quasi-geostrophic equations.
The effect of the CE on 3D complex flows also needs to be studied.
Finally, we plan to investigate whether the CE has a significant effect on spatially-filtered ROMs that are different from the DDC-ROM considered in this paper, e.g., the physically constrained data-driven ROM~\cite{mohebujjaman2018physically} or the approximate deconvolution ROM~\cite{xie2017approximate}.

%


\bibliographystyle{spmpsci}      
\bibliography{traian}   

%
%

\end{document}